\documentstyle[12pt]{article} 
\begin{document}
\newcommand{\singlespace}{
    \renewcommand{\baselinestretch}{1}
\large\normalsize}
\newcommand{\doublespace}{
   \renewcommand{\baselinestretch}{1.2}
   \large\normalsize}
\renewcommand{\theequation}{\thesection.\arabic{equation}}

\input amssym.def
\input amssym
\setcounter{equation}{0}
\def \ten#1{_{{}_{\scriptstyle#1}}}
\def \Z{\Bbb Z}
\def \C{\Bbb C}
\def \R{\Bbb R}
\def \Q{\Bbb Q}
\def \N{\Bbb N}
\def \l{\lambda}
\def \V{V^{\natural}}
\def \wt{{\rm wt}}
\def \tr{{\rm Tr}}
\def \Res{{\rm Res}}
\def \End{{\rm End}}
\def \Aut{{\rm Aut}}
\def \mod{{\rm mod}}
\def \Hom{{\rm Hom}}
\def \im{{\rm im}}
\def \<{\langle} 
\def \>{\rangle} 
\def \w{\omega}
\def \o{\omega}
\def \t{\tau }
\def \ch{{\rm ch}}
\def \a{\alpha }
\def \b{\beta}
\def \e{\epsilon }
\def \la{\lambda }
\def \om{\omega }
\def \O{\Omega}
\def \qed{\mbox{ $\square$}}
\def \pf{\noindent {\bf Proof: \,}}
\def \voa{vertex operator algebra\ }
\def \voas{vertex operator algebras\ }
\def \p{\partial}
\def \1{{\bf 1}}
\def \g{{\frak g}}
\singlespace
\newtheorem{thmm}{Theorem}
\newtheorem{th}{Theorem}[section]
\newtheorem{prop}[th]{Proposition}
\newtheorem{coro}[th]{Corollary}
\newtheorem{lem}[th]{Lemma}
\newtheorem{rem}[th]{Remark}
\newtheorem{de}[th]{Definition}
\newtheorem{con}[th]{Conjecture}
\newtheorem{ex}[th]{Example}

\begin{center}
{\Large {\bf Holomorphic Vertex Operator Algebras of Small Central Charges}} \\
\vspace{0.5cm}

Chongying Dong\footnote{Supported by NSF grant 
 DMS-9987656 and faculty research funds granted by the University of 
California at 
Santa Cruz (dong@math.ucsc.edu).} and 
Geoffrey Mason\footnote{Supported by NSF grant DMS-9700909
and faculty research funds granted by the University of 
California at 
Santa Cruz (gem@cats.ucsc.edu).}
\\
Department of Mathematics, University of
California, Santa Cruz, CA 95064
\end{center}
\hspace{1.5 cm}

\begin{abstract} 
We provide a rigorous mathematical foundation to the study 
of strongly rational, holomorphic vertex operator algebras $V$ of 
central charge $c = 8, 16$ and $24$ initiated by Schellekens. If $c = 8$ 
or $16$ we show that $V$ is isomorphic to a lattice theory 
corresponding to a rank $c$ even, self-dual lattice. If $c = 24$ we 
prove, among other things, that either $V$ is isomorphic to a lattice 
theory corresponding to a Niemeier lattice or the Leech lattice, or 
else the Lie algebra on the weight one subspace $V_1$ is semisimple 
(possibly $0$) of Lie rank less than $24$.
\end{abstract}

\section{Introduction}

One of the highlights of discrete mathematics is the
classification of the positive-definite, even, unimodular lattices $L$
of rank at most $24$, due originally to Minkowski, Witt and Niemeier
(cf. [CS] for more information and an extensive list of references).
One knows (loc. cit.) that such a lattice has rank divisible by $8$;
that the $E_8$ root lattice is (up to isometry) the unique example of
rank $8$; that the two lattices of type $E_8 + E_8$ and $\Gamma_{16}$ are
the unique examples of rank $16$; and that there are $24$ inequivalent
such lattices of rank $24$. In each case, the lattice may be
characterized by the nature of the semi-simple root system naturally
carried by the set of minimal vectors (i.e. those of squared length
$2$). 

The theory of vertex operator algebras is a newer subject which enjoys
several parallels with lattice theory. Already in [B], Borcherds
pointed out that one can naturally associate a vertex operator algebra
$V_L$ to any positive-definite, even lattice $L$ (cf. [FLM] for a
complete discussion). It is known that $V_L$ is rational and that it
is holomorphic precisely when $L$ is self-dual ([D], [DLM1]) (we defer
the formal introduction of technical definitions concerning vertex
operator algebras until Section 2). Since the central charge $c$ of
the vertex operator algebra $V_L$ is precisely the rank of $L$, the
classification of holomorphic vertex operator algebras of central
charge at most $24$ may be construed as a generalization of the
corresponding problem for even, unimodular lattices. We will say that
a holomorphic vertex operator algebra has {\em small central charge}
in case it satisfies $c\leq 24.$

Schellekens was the first to consider the problem of classifying 
holomorphic vertex operator algebras $V$ of small central 
charge [Sch]. Based on extensive computation, Schellekens wrote down a 
list of $71$ integral $q$-expansions $q^{-1} + constant + 196884q + 
...$ and, among other things, conjectured that the graded character 
of a $V$ satisfying $c=24$ is necessarily equal to one of these $71$ 
$q$-expansions. As is well-known,  it is only the constant 
term that distinguishes the $q$-expansions from each other, the 
constant in question being the dimension of the Lie algebra naturally 
defined on the weight one subspace $V_1$ of $V$ (see below for more 
details). Schellekens in fact wrote down a list of $71$ Lie algebras 
(including levels) which are the candidates for $V_1$. (It turns out 
that if $V$ has central charge strictly less that $24$ then the 
graded character of $V$ is uniquely determined, so at least as far as 
the dimension of the weight one subspace is concerned, the case 
$c=24$ carries the most interest.)

The purpose of the present paper is to put the Schellekens program on a
firm mathematical foundation within the general context of rational
conformal field theory, and to make a start towards the classification
(up to isomorphism) of the holomorphic vertex operator algebras of
small central charge. In addition to unitarity, there are other
(unstated) assumptions in [Sch]. This circumstance means that we cannot
assume the results of (loc. cit.) and need to find new approaches. In
some ways we will go much further than [Sch] in that we will be able to
give an adequate characterization of the holomorphic lattice theories
among all holomorphic theories of small central charge. On the other
hand, although we are able to establish numerical restrictions on the
nature of the Lie algebra on $V_1$ which show that there are only a
finite number of possibilities, we are at present unable to show, assuming 
$c=24,$ that it is necessarily one of the $71$ on Schellekens' list. 
The remaining
obstacle is essentially that of establishing that the levels of the
associated Kac-Moody Lie algebras are positive integers (which is
immediate if unitarity is assumed).

We now state some of our main results, and for this purpose we will
take $V$ to be a holomorphic vertex operator algebra of CFT type (see
Section 2) which is $C_2$-cofinite. (In the language of [DM], $V$ is
strongly rational and holomorphic.). By results of Zhu [Z], this
implies that $c$ is a positive integer divisible by $8$.

\begin{thmm}\label{t1} Suppose that $c=8$. Then $V$ is isomorphic to the 
lattice theory $V_{E_8}$ associated to the $E_8$ root lattice.
\end{thmm}

\begin{thmm}\label{t2} Suppose that $c=16$. Then $V$ is isomorphic to a lattice
theory $V_L$ where $L$ is one of the two unimodular rank $16$
lattices.
\end{thmm}

Thus for the cases $c=8$ and $16$, the classification of the holomorphic
vertex operator algebras completely mirrors that of the corresponding
lattices. These two theorems are commonly assumed in the physics
literature, and are related to the uniqueness of the heterotic string
([Sch], [GSW]).

\begin{thmm}\label{t3} Suppose that $c=24$. Then the Lie algebra on $V_1$ is
reductive, and exactly one of the following holds:

(a) $V_1 = 0$.

(b) $V_1$ is abelian of rank $24$. In this case $V$ is isomorphic to the
lattice theory $V_{\Lambda}$ where $\Lambda$ is the Leech lattice.

(c) $V_1$ is a semi-simple Lie algebra of rank $24$. In this case $V$ is
isomorphic to the lattice theory $V_L$ where $L$ is the even,
unimodular rank $24$ lattice whose root system is the same as the root
system associated to $V_1$.

(d) $V_1$ is a semi-simple Lie algebra of rank less than $24$.
\end{thmm}

We actually establish more than is stated
here. It is a well-known conjecture [FLM] that the Moonshine module
is characterized among all $c=24$ holomorphic theories by the condition
$V_1 = 0$. Our methods are less effective when there is no Lie algebra,
however we will show that in this case $V_2$ carries the structure of a
simple, commutative algebra (of dimension 196,884). The commutativity
and dimension formula are well-known; it is the simplicity that is
novel here. (The inherent difficulty in dealing with $V_2$ is that it
is not an associative algebra, indeed it is not even power associative,
and there seem to be no useful identities which are satisfied.) In case
(d) we show that the simple components ${\frak g}_i$ of $V_1$ have levels $k_i$
and dual Coxeter numbers $h_i^{\vee}$ such that the identity
\begin{equation}\label{1.1}
\frac{h_i^{\vee}}{k_i}  = \frac{ \dim V_1 -24}{24}
\end{equation}                                            
holds for each $\g_i$. This implies that there are only finitely many choices for the
family of pairs $(\g_i, k_i)$ determined by $V_1$. Note that the
condition (\ref{1.1}) was already identified by Schellekens [Sch].  
We are in fact able to extract some further numerical 
restrictions on the levels $k_i$, but these fall well short of the 
expectation that they are all positive integers, and we forgo any 
discussion of this beyond (\ref{1.1}). We also establish that if $V_1$ is 
semisimple in Theorem \ref{t3}, then the Virasoro element in $V$ coincides 
with the usual Virasoro element associated with an affine Lie 
algebra defined by the Sugawara construction. We expect the result to be useful in 
further analysis of the situation.

The main inspiration for the proof of our results originates from our 
recent paper [DM]. There, we introduced methods based on the theory 
of modular forms used in tandem with techniques from vertex operator 
algebra theory. In particular we obtained a simple numerical 
characterization of the lattice vertex operator algebras among all 
rational vertex operator algebras. The present paper is in many ways 
a continuation and elaboration of [loc.cit.]. When the central 
charge is small, the theory of modular forms gives very precise 
numerical information about the vertex operator algebra which allows 
us to show, under certain circumstances, that the numerical 
characterizations of [DM] are applicable. This leads to Theorems 
1-3. Modular-invariance also underlies the other results that we 
obtain.

In addition to the Schellekens program that we have already 
discussed, another potential application of Theorem \ref{t3} is to the FLM 
conjecture regarding the Moonshine module alluded to above. Namely, 
Theorem \ref{t3} shows that the Leech lattice theory $V_{\Lambda}$ is the only 
$c=24$ holomorphic theory for which the Lie algebra $V_1$ is both 
non-zero and not semi-simple. On the other hand, the Moonshine module 
is closely related to $V_{\Lambda}$, being built from it by a 
$\Z_2$-orbifold construction [FLM].

The paper is organized as follows: we gather together some
preliminaries in Section 2, and prove Theorems \ref{t1} - \ref{t3}
together with the supplementary result in case(a) of Theorem 3 in
Section 3. In section 4 we identify the Virasoro element with the
Sugawara construction.

\section{Preliminary results}
\setcounter{equation}{0}

For general background on the theory  of vertex operator algebras, we refer 
the reader to [FLM], [FHL]. As usual, for a state $v$ in a vertex operator 
algebra $V$, we 
denote the corresponding vertex operator by
$$Y(v, z)  =  \sum_{n \in \Z}v_n z^{-n-1},$$
while the vertex operator corresponding to the conformal (Virasoro) vector is
$$Y(\omega, z)  =  \sum_{n \in \Z} L(n)z^{-n -2}.$$

$V$ is called {\em rational} if all admissible $V$-modules are completely 
reducible (cf. [DLM1] for the definition of admissible module). It 
was shown in [DLM2] that this implies that $V$ has only finitely many 
inequivalent simple modules. A rational vertex operator algebra is 
called {\em holomorphic}  if it has a unique simple module, namely the 
adjoint module $V$. $V$ is said to be $C_2$-{\em cofinite} in case the 
subspace spanned by elements $u_{-2}v$ for $u, v \in V$ is of finite 
codimension. Finally, we say that $V$ is of {\em CFT-type} in case the 
natural $\Z$-grading on $V$ takes the form
\begin{equation}\label{2.1}
V  =  V_0 \oplus V_1  \oplus \cdots\ \   {\rm with}\ \   V_0  = \C\1 .
\end{equation}

Throughout the rest of this paper, we assume that $V$ is a 
$C_2$-cofinite, holomorphic vertex operator algebra of CFT-type of 
small central charge $c\leq 24$. In order to avoid the case of the trivial 
vertex operator algebra $V = \C$, we also assume that $V$ has 
dimension greater than one.

Next we discuss some consequences of these assumptions for the 
structure of $V$. Many of them are well-known. First, the weight $1$ 
subspace $V_1$ of $V$ carries a natural structure of Lie algebra 
given by
\begin{equation}\label{2.2}
[u,v]  =  u_0v
\end{equation}
for $u,v \in V_1$.   
Because the adjoint module is the unique simple $V$-module, then the 
contragredient module $V'$ is necessarily isomorphic to $V$. This is 
equivalent to the existence of a non-degenerate, invariant bilinear 
pairing
$$ \<  ,  \>:  V \times V \to \C$$
which is necessarily symmetric. (For the theory of contragredient 
modules, cf. Section 3 of [FHL].) Because of (\ref{2.1}), Li's theory of 
invariant bilinear forms [L1] shows that we have $L(1)V_1 = 0$ and 
that $\<  ,  \>$ is uniquely determined up to an overall scalar. It is 
convenient to fix the normalization so that
\begin{equation}\label{2.3}
\< \1,\1 \>  =  -1.
\end{equation}
In particular, the restriction of $\<  ,  \>$ to $V_1$ is given by
\begin{equation}\label{2.4}
\<u, v\>\1  =  u_1v
\end{equation}
for $u, v \in V_1.$ 

In the language of [DM], it follows from what we have said that $V$ 
is {\em strongly rational,} so that the results of (loc. cit.) apply. They 
tell us that the following hold:

(I) $V_1$ is  a reductive Lie algebra of Lie rank $l\leq c.$ 

(II) $l = c$ if, and only if, $V$ is isomorphic to a lattice 
theory $V_L$ for some   positive-definite, even, unimodular lattice $L$. 

We also note that by results of Zhu [Z] (cf. [DLM3]), $c$ is 
necessarily a positive integer divisible by $8$. 
So in fact $c = 8, 16$, or $24.$

We refer the reader to [Z] and [DM] for an extended discussion of 
the role of modular-invariance in the theory of rational vertex 
operator algebras. We need  to recall the genus
one vertex operator algebra $(V,Y[\ ],\1, \omega-c/24)$ from [Z]. 
The new vertex 
operator associated to  a homogeneous element $a$ is given by
\[ Y[a,z] = \sum_{n\in\Z}a[n]z^{-n-1} = Y(a,e^{z} -1)e^{z\wt{a}}
\] while the Virasoro element is $\tilde{\omega} = \omega-c/24$.
Thus 
\[
a[m] = \Res_z\left(Y(a,z)(\ln{(1+z)})^m(1+z)^{\wt{a}-1}\right)
\]
and \[
a[m] = \sum_{i=m}^\infty c(\wt{a},i,m)a(i)
\]
for some scalars $c(\wt{a},i,m)$ such that $c(\wt{a},m,m)=1.$ 
In particular,
\[
a[0]=\sum_{i\geq 0}{\wt{a}-1 \choose i}a(i).
\]
We also write
\[
L[z] = Y[\omega-c/24,z] = \sum_{n\in\Z} L[n]z^{-n-2}.
\]
Then the $L[n]$ again generate a copy of the Virasoro algebra with
the same central charge $c.$ Now $V$ is graded by
the $L[0]$-eigenvalues, that is
\[
V=V_{[0]}\oplus V_{[1]}\oplus \cdots 
\]
where $V_{[n]}=\{v\in V|L[0]v=nv\}.$ 

We will need the following facts: Let $v \in V$ 
satisfy $L[0]v = kv.$ 
Then the graded trace
\begin{equation}\label{2.7}
Z(v, \tau)  =  \tr_V o(v) q^{L(0) - c/24}  =  q^{-c/24}\sum_{n\geq 0} 
\tr_{V_n}o(v) q^n 
\end{equation}
is a modular form on $SL(2,\Z)$ (possibly with character), holomorphic 
in the complex upper half-plane $H$ and of weight $k$. Here, we have 
set $o(v)$ to be the {\em zero mode} of $v$, defined via $o(v) = v_{\wt v - 
1}$ if $v$ is homogeneous, and extended by linearity to the whole of 
$V$. Moreover, $\tau$ will denote an element in $H$ and $q = e^{2 \pi 
i \tau}$. In particular, if we take $v$ to be the vacuum element then 
(\ref{2.7}) is just the graded trace
\begin{equation}\label{2.8}
\ch_qV=q^{-c/24} \sum_{n\geq 0} (\dim V_n) q^n
\end{equation}
and is a modular function of weight zero on $SL(2,\Z)$. Because of the 
holomorphy of $\ch_q(V)$ in $H$ and our assumption that $c\leq 24$, one 
knows ([K], [Sch]) that (\ref{2.8}) is uniquely determined up to an additive 
constant, and indeed is determined uniquely if $c\leq 16$. The upshot is 
this:

\begin{lem}\label{l2.1} One of the following holds:

(a) $c = 8$ and $\ch_q(V)  = \Theta_{E_8}(q)/ \eta(q)^8 =  q^{-1/3}(1 + 248q +
 \cdots)$

(b) $c = 16$ and $\ch_q(V)  =  (\Theta_{E_8}(q))^2/ \eta(q)^{16}  = 
q^{-2/3}(1 + 496q + \cdots)$

(c) $c = 24$ and $\ch_q(V) = J(q) + const = q^{-1} + const + 196884q +\cdots.$
\end{lem} 

Here, we have introduced the theta function $\Theta_{E_8}(q)$ of 
the $E_8$ root lattice
$$\Theta_{E_8}(q)=\sum_{\alpha\in E_8}q^{(\alpha,\alpha)/2}$$
(where $E_8$ denotes the root lattice of type $E_8$ normalized so that
the squared length of a root is $2$) 
 as well as the eta function 
$$\eta(q)=q^{1/24}\prod_{n\geq 1}(1-q^n);$$
$$J(q) = q^{-1} + 0 + 196884q +\cdots $$
is the absolute modular invariant normalized to 
have constant term zero (alias the graded character of the Moonshine 
Module [FLM]). We also need the 'unmodular' Eisenstein series of 
weight two, namely
$$E_2(q)  =  1 - 24 \sum_{n\geq 1}\sigma_1(n) q^n$$
where $\sigma_1(n)$ is the sum of the divisors of $n.$

\begin{lem}\label{l2.2} For states $u, v$ in $V_1$ we have
          $$Tr_V o(u)o(v) q^{L(0) - c/24}  =  \frac{   \< u, v\>}{24} 
(\frac{ 24}{c}D_q(ch_qV)  + E_2(q)ch_qV)$$
where $D_q = q\frac{ d}{ dq}$.
\end{lem}

Proof: First recall the important identity, Proposition 4.3.5 in [Z]. 
If we pick a pair of states
$u, v \in V_1$ this result yields the following:
\begin{equation}\label{2.9}
\tr_V o(u)o(v) q^{L(0) - c/24}  =  Z(u[-1]v, \tau)  + 1/ 12 
E_2(\tau) Z(u[1]v, \tau).
\end{equation}
The term $Z(u[-1]v, \tau)$ is a modular form of weight 2, and we shall 
be able to write it down explicitly. Indeed, the l.h.s. of (\ref{2.9}) has 
leading term $\kappa(u,v)q^{1-c/ 24}$ where
$\kappa(u,v)$ denotes the usual Killing form on $V_1$, whereas the 
leading term of the second summand on the r.h.s. of (\ref{2.9}) is equal to 
$\frac{-1}{12}\<u, v\>q^{-c/ 24}$. From this discussion we conclude that the 
leading term of $Z(u[-1]v, \tau)$ is equal to
$\frac{1}{12}\<u,v\>q^{-c/24}$. But up to scalars, the unique form of weight 
$2$ on $SL(2,\Z)$ which is holomorphic in $H$ and has a pole of order 
$c/ 24$ at infinity is $D_q(ch_q(V))$. As a result, we see that the 
first term on the r.h.s. of (\ref{2.9}) is equal to $\frac{2\<u, v\>}{c} 
D_q(ch_q(V)).$ The lemma follows from this, (\ref{2.9}), and (\ref{2.4}).
\qed

Identifying coefficients of $q^{1 - c/24}$ in the formula of Lemma
\ref{l2.2} yields

\begin{coro}\label{c2.3} $\kappa(u,v)  = 2\<u , v\> (\frac{\dim V_1}{c} - 1)$.
\end{coro}

\section{Proof of the Theorems}
\setcounter{equation}{0}

 We consider the three possibilities 
$c = 8, 16, 24$ in turn. In the first two cases the idea is to show 
that (II) {\em always} applies, while in the third case we study {\em when} it 
applies.

Case 1: $c = 8$. In this case dim $V_1 = 248$ by Lemma \ref{l2.1} (a), so 
that $\kappa(u, v) = 60\<u, v\>$ by Corollary \ref{c2.3}. Since $\< , \>$ is 
non-degenerate then so too is the Killing form. We conclude that in 
fact $V_1$ is semi-simple of dimension $248$, and by (I) the Lie 
rank is no greater than $8$. By the classification of semi-simple Lie 
algebras, we conclude that in fact $V_1$ is the Lie algebra of type 
$E_8$ and Lie rank $8$. Now (II) and the fact that there is a unique 
positive-definite, even,  unimodular lattice of rank $8$, namely the 
$E_8$ root lattice, completes the proof of Theorem \ref{t1}.

Case 2: $c = 16$. This is very similar to case 1. Namely, we have dim 
$V_1 = 496$ by Lemma \ref{l2.1} (b), whence $\kappa(u,v) = 60\<u, v\>$ once 
more. So again $V_1$ is semi-simple, and we proceed as in case 1, now 
using the fact that there are just the two positive-definite, even, 
unimodular lattices. Theorem \ref{t2} follows. 

Case 3: $c = 24$. Set $dim V_1 = d$. From Corollary 2.3 we see in 
this case that
\begin{equation}\label{3.1}
\kappa(u, v) = \<u, v\> (d - 24)/12.
\end{equation}

If $d = 24$  then the Killing form is identically zero. Then $V_1$ is 
solvable by Cartan's criterion and therefore abelian since $V_1$ is 
in any case reductive by (I). So if $d = 24$ then we have shown 
that $V_1$ has rank $24$ and is therefore a lattice theory $V_L$ for 
suitable $L$ by (II). The fact that $(V_L)_1$ is abelian tells us 
that $L$ has no roots i.e., no vectors of squared length two, and 
that $L$ is therefore the Leech lattice (cf. [CS]). This deals with 
case (b) of Theorem \ref{t3}.

Let us assume from now on that $d \ne 0,24$. Together with (\ref{3.1}), 
this tells us that the Killing form on $V_1$ is again non-degenerate, 
so that $V_1$ is a semi-simple Lie algebra of Lie rank $l$ no greater 
than $24$ by (I). Moreover, if $l = 24$ then $V$ is a lattice 
theory by (II) once more. This confirms part (c) of Theorem \ref{t3}.
\qed

Next we consider the levels of the affine Lie algebras spanned by the 
vertex operators $Y(u, z)$, $u$ in $V_1$. For states $u, v \in V_1$ 
and integers $m, n$ we have
\begin{equation}\label{a1}
[u_m, v_n]  =  (u_0v)_{m + n}  +  m u_1v \delta_{m, - n},
\end{equation}
whereas the usual relations for a Kac-Moody Lie algebra of level $k$ 
associated to a simple Lie algebra $\g$ take the form
\begin{equation}\label{a2}
[a_m, b_n] = [a, b]_{m + n} + k(a, b)m\delta_{m, - n}
\end{equation}
where $(a, b)$ is the non-degenerate form on $\g$ normalized so that 
$(\alpha, \alpha) = 2$ for a long root $\alpha$. 

Let $V_1$ be a direct sum 
\begin{equation}\label{6.11}
V_1  =  {\frak g}_{1, k_1}  \oplus {\frak g}_{2, k_2} \oplus ... \oplus 
{\frak g}_{n,k_n}
\end{equation}
of simple Lie algebras ${\frak g}_i$ whose corresponding affine Lie
algebra has level $k_i.$ By comparing (\ref{a1}) and (\ref{a2}), using
Corollary \ref{c2.3} we obtain for $u,v\in{\frak g}_{i}$ that
\begin{equation}\label{3.2} 
 \kappa_{{\frak g}_{i}}(u, v) =( d-24) k_i(u, v)/ 12
\end{equation}
where $\kappa_{{\frak g}_{i}}$ denotes the restriction of the Killing 
form to ${\frak g}_{i}.$ Now $\kappa_{{\frak g}_{i}}(h_{\alpha}, h_{\alpha}) = 4 h_i^{\vee} $ for a long root $\alpha$, where 
$h_i^{\vee}$ is the dual Coxeter number of the root system associated to ${\frak g}_{i}$. Therefore, (\ref{3.2}) tells us that for each simple 
component ${\frak g}_{i}$ of 
$V_1$, of level $k_i$ and dual Coxeter number $h_i^{\vee}$, the ratio
\begin{equation}\label{3.3}
h_i^{\vee}/ k_i  =(d- 24)/ 24
\end{equation}
is independent of ${\frak g}_{i}$.

\begin{prop}\label{p3.1} Assume that case (a) of Theorem \ref{t3} holds. Then $B = 
V_2$ carries the structure of a (non-associative) simple, commutative 
algebra with respect to the product $a.b = a_1b$.
\end{prop}

\pf We take for granted the well-known facts that $B$ is indeed a 
non-associative, commutative algebra with respect to the indicated 
product, and that the pairing $\< ,  \>: B \times B \to \C$ defined by
\begin{equation}\label{3.4} 
a_3b = \<a,b\>\1
\end{equation}
endows $B$ with a non-degenerate, invariant trace form. That is, $\< 
,  \>$ is symmetric and satisfies $\<ab, c\> = \<a , bc\>$. Moreover $B$ 
has an identity element $1/2\omega$. Set $d = \dim B = 196884$.

Next we state two more results that we will need. Each can be 
established using modular-invariance arguments along the same lines 
as before. Alternatively, we may use results of Section 4 of [M]:
\begin{equation}\label{3.5}
\tr_B o(ab)  =(d/3) \<a, b\>.
\end{equation}
If $e^2 = e$ is in $B$ then 
\begin{equation}\label{3.6}
\tr_B o(e)^2 = 4620 \<e, e\> + 20336 \<e,e\>^2.
\end{equation}

Turning to the proof of the Proposition, we first show that $B$ 
is semisimple. Indeed, (\ref{3.5}) guarantees that the form $\tr_B o(ab)$ is 
non-degenerate, and this is sufficient to establish that $B$ is 
semi-simple. To see this, recall a well-known result of Dieudonne (cf.
[S, Theorem 2.6]) that an arbitrary algebra $B$ is semisimple if it 
has a non-degenerate trace form and contains no non-zero nilpotent 
ideals. We will show that indeed there are no non-zero nilpotent 
ideals in $B$. If not, we may choose a minimal non-zero nilpotent 
ideal $M$ in $B$. Note that $M^2 = 0$. Let $m$ be a non-zero element 
in $M$, and let $b$ in $B$ be arbitrary. Then $mb$ lies in $M$, so 
that $(mb)B\subset M$ and $(mb)M = 0$. This shows that each element $mb$ 
is nilpotent as a multiplication operator on $B$, and hence $\tr_B mb 
= 0$ for all $b$. This contradicts the non-degeneracy of $\< , \>$.

From the last paragraph we know that $B$ can be written as an 
(orthogonal) direct sum of simple ideals
 $$B  =  B_1 + B_2 + \cdots + B_t.$$
We must show that $t = 1$. Write $1/2 \omega = e_1 + e_2 + ... + e_t$ 
where $e_i$ is the identity element of $B_i$. In particular, $e_i$ is 
a central idempotent of $B$. By (\ref{3.5}) and (\ref{3.6}) we obtain
\begin{equation}\label{3.7}
(d/3)\<e, e\> = 4620\<e, e\> +20336\<e, e\>^2
\end{equation}
where $e$ is any of the idempotents $e_i$. Notice that (\ref{3.7}) yields 
that
\begin{equation}\label{3.8}
\<e, e\>= 3. 
\end{equation}

It is well-known that (\ref{3.8}) implies that the components of the vertex 
operator $Y(2e, z)$ generate a Virasoro algebra of central charge 
$24$, so that the total central charge must be $24t$. Hence, $t = 1$ 
as required. \qed

\section{Some Virasoro elements}
\setcounter{equation}{0}

Throughout this section we assume that $V$ is a strongly rational, 
holomorphic vertex operator algebra of central charge $24$ as before, 
and we assume in addition that $V_1$ is a (non-zero) semisimple Lie 
algebra of Lie rank $l.$ Let $V_1$ have decomposition (\ref{6.11}) 
into simple Lie algebras. Set $d_i = \dim {\frak g}_{i,k_i}.$ Let $( , )$ 
denote the normalized 
invariant bilinear form on ${\frak g}_{i,k_i}$ with the property that 
$(\alpha,\alpha)  =  2$ for a long root $\alpha$ in ${\frak g}_{i,k_i}.$  
It is known (cf. [DL], [FZ], [K], [L2]) that the element
\begin{equation}\label{4.1}
\omega_i=\frac{1}{2(k_i+h_i^{\vee})}\sum_{j=1}^{d_i}u^j_{-1}u^j_{-1}{\bf 1}
\end{equation}
is a Virasoro element of central charge $c_i=\frac{k_i\dim {\frak
g}_i}{k_i+h_i^{\vee}},$ where $u^1, u^2, ... u^{d_i}$ is an
orthonormal basis of ${\frak g}_i$ with respect to $( , )$.

There are three Virasoro elements in $V$ that are relevant to a 
further analysis of the situation. Namely, in addition to the 
original Virasoro element $\omega$ in $V$, we have
\begin{equation}\label{4.2}
\omega_{aff}  = \sum_{i=1}^n\omega_i
\end{equation}
and
\begin{equation}\label{4.3}
\omega_H  = \frac{1}{2}\sum_{i=1}^l(h^i_{-1})^{2}{\bf 1} 
\end{equation}
where $h^{1}, ... , h^{l}$ is an orthonormal basis of a maximal
abelian subalgebra $H$ of $V_1$ with respect to the inner product $\<,
\>.$ Note that as a consequence of equation (\ref{3.3}),
$\omega_{aff}$ has central charge $\sum c_i = 24$. We omit further
discussion of $\omega_H$, but will prove:
\begin{prop}\label{p4.1} $\omega_{aff}  =  \omega$.
\end{prop}

\pf  Consequences of modular-invariance again underlie the proof of 
the Proposition, notably the absence of cusp-forms of small weight on 
$SL(2,\Z)$. We introduce some notation for Virasoro operators: in 
addition to the usual operators $L(n)$ associated to $\omega$ , we use 
$L^{aff}(n)$ for the corresponding operators associated to $\omega_{aff}.$ 
We also set $\omega' = \omega - \omega^{aff}$. We will soon see that 
$\omega'$  is itself a Virasoro element, and define its component 
operators to be  $L'(n)$. We eventually want to prove of course that 
$\omega' = 0$. We proceed in a series of steps.

Step 1: $\omega'$  is a highest weight vector of weight 2 for the 
Virasoro algebra generated by ${L(n)}$.

The definition of $\omega'$  shows that it lies in $V_2$, so it 
suffices to show that $\omega'$  is annihilated by the operators 
$L(1)$  and $L(2)$. First calculate that if $u$ is an element of 
$V_1$  then $[L(1), u_{-1}]  =  u_0$. Then from the definitions it 
follows easily that $L(1)$ annihilates both $\omega_{aff}$  and 
$\omega_H$, whence it also annihilates $\omega'$.  To establish that 
$L(2)$ annihilates $\omega'$ we must show that
\begin{equation}\label{4.4}
L(2)\omega_{aff}  =  L(2)\omega_H =12.
\end{equation}
Now for  $u$  in $V_1$  we have $[L(2), u_{-1}]  =  u_1$. Then
$$L(2) u^{j}_{-1} u^{j}_{-1}{\1} =  - \<u^{j}, u^{j}\>  =k_i,$$
and (\ref{4.4}) follows.

Step 2: $\omega'$ is a Virasoro vector of central charge $0$.

Use step 1, in particular that $L(1)$ annihilates $\omega_{aff}$, 
to see that
$$[L(m), L^{aff}(n)]  =  (m-n) L^{aff}(m+n)  -2(m^{3} 
- m)\delta_{m, -n}.$$
The result follows easily from this.

Step 3: $Z_V(\omega', \tau)  =  0$.

Since $\omega'$  is a highest weight vector of weight 2 for the 
Virasoro operators corresponding to $\omega$ then from our earlier 
discussion of (\ref{2.7}) we see that
$Z_V(\omega', \tau)$ is a modular form on $SL(2,Z)$ of weight 2 and is 
holomorphic in the upper half-plane. Moreover there is no pole in the 
q-expansion of $Z_V(\omega', \tau)$, so $Z_V(\omega', \tau)$  is in fact 
a holomorphic modular form of weight 2 on $SL(2,Z)$, hence must be 
zero.

Step 4: $\tr_VL'(0)^2q^{L(0)-1}=0.$

By Proposition 4.3.5 in [Z] we have 
$$\tr_VL'(0)^2q^{L(0)-1}=Z_V(L'[-2]\omega',\tau)-\sum_{k\geq
1}E_{2k}(\tau) Z_V(L'[2k-2]\omega',\tau)$$ where
$Y[\omega',z]=\sum_{m\in\Z}L'[m]z^{-m-2}$ and the functions
$E_{2k}(\tau)$ are Eisenstein series of weight $2k,$ normalized as in
[DLM3].  $E_{2k}(\tau)$ is a holomorphic modular form on
$SL(2,\Z)$ if
$k>1$. Since $L'[2k-2]\omega'=0$
if $k\geq 2,$ $L'[0]\omega'=2\omega',$ and $L[0]L'[-2]\omega'=4L'[-2]\omega',$
we see that
$$\tr_VL'(0)^2q^{L(0)-1}=Z_V(L'[-2]\omega',\tau)$$ 
is a modular form of weight 4 for $SL(2,\Z).$  Since $L'(0)V_1=0,$ 
it is in fact a cusp form, hence equal to zero.

Step 5. All the eigenvalues of $L^{aff}(0)$ on $V$ are real.

In order to see this recall the decomposition (\ref{6.11}) and set $\g=V_1.$ 
Consider the
affine Lie algebra
$$\hat{\frak g}={\frak g}\otimes \C[t,t^{-1}]\oplus \C$$
with bracket
$$[x\otimes t^p,y\otimes t^q]=[x,y]\otimes t^{p+q}+p\delta_{p+q,0}\<x,y\>$$
for $x,y\in {\frak g}$ and $p,q\in\Z.$ Since each $V_m$ is 
finite dimensional we see that $V$ has a composition series 
as a module for $\hat{\frak g}$ such that each factor
is an irreducible highest weight $\hat{\frak g}$-module. So it is enough to 
show that $L^{aff}(0)$ has only real eigenvalues
on any irreducible highest weight $\hat{\frak g}$-module. Note that
such an irreducible  highest weight $\hat{\frak g}$-module is a tensor
product of irreducible highest weight $\hat{\frak g}_i$-modules
$L(k_i,\Lambda_i)$ of level $k_i$ ($i=1,...,n$) for some
dominant weight $\Lambda_i$ in the weight lattice
of ${\frak g}_i$ as $V$ is a completely reducible  ${\frak g}_i$-module.
Here $L(k_i,\Lambda_i)$ is the unique irreducible quotient 
of the generalized Verma module $U(\hat{\frak g}_i)\otimes_{U({\frak g}_i\otimes \C[t]+\C)}L(\Lambda_i)$ where $L(\Lambda_i)$ is 
the highest weight module for ${\frak g}_i$ with highest 
weight $\Lambda_i$ and $x\otimes t^m$ acts as zero if $m>0$ and
$x\otimes t^0$ acts as $x.$ So it is enough
to show that $L_i(0)$ has only real eigenvalues on $L(k_i,\Lambda_i)$
where $Y(\omega_i,z)=\sum_{m\in\Z}L_i(m)z^{-m-2}.$ 
   
It is well-known that the eigenvalues of $L_i(0)$ on $L(k_i,\Lambda_i)$
are the numbers 
$\frac{(\Lambda_i+2\rho_i,\Lambda)}{2(k_i+h_i^{\vee})}+m$ (cf. [DL],
[K]) for nonnegative integers $m$ where $\rho_i$ is the half-sum of the
positive roots of ${\frak g}_i.$ Since $k_i$ is rational, it is clear that
$\frac{(\Lambda_i+2\rho_i,\Lambda)}{2(k_i+h_i^{\vee})}+m$ is real, as required.

Step 6. $\omega'=0.$

By Step 5, all eigenvalues of $L'(0)$ are real on $V.$ This implies
that $L'(0)^2$ has only nonnegative eigenvalues. By Step 4  we conclude that all the
eigenvalues of $L'(0)$ are zero.  Since $L'(0)\omega'=2\omega',$
we must have $\omega'=0,$ as desired. \qed

\end{document}